\documentclass[12pt]{amsart}
\usepackage{amsxtra}
\usepackage{mathrsfs}

\newtheorem{theorem}{Theorem}
\newtheorem{definition}{Definition}

\renewcommand{\ni}{\noindent}

\newcommand{\ZZ}{\ensuremath{\mathbb{Z}}}
\newcommand{\QQ}{\ensuremath{\mathbb{Q}}}

\newcommand{\CC}{\ensuremath{\mathbb{C}}}
\newcommand{\FF}{\ensuremath{\mathbb{F}}}

\newcommand{\tensor}{\otimes}

\newcommand{\ov}{\overline}
\newcommand{\eps}{\varepsilon}

\newcommand{\GL}{\operatorname{GL}}
\newcommand{\Gal}{\operatorname{Gal}}

\title{Rigid Calabi-Yau Threefolds over \QQ\ Are Modular}

\author{{\sc Fernando Q. Gouv\^ea} and {\sc Noriko Yui}}
\thanks{This work was partially supported by
        the Natural Sciences and Engineering Research Council
        of Canada (NSERC)}
\date{}

\begin{document}

\begin{abstract}
The proof of Serre's conjecture on Galois representations over finite
fields allows us to show, using a method due to Serre himself, that all
rigid Calabi-Yau threefolds defined over \QQ\ are modular.
\end{abstract}

\maketitle



\vspace{2\baselineskip}

In the mid-1980s, J.-P. Serre conjectured in \cite{duke} that all
absolutely irreducible odd two-dimensional representations of
$G_\QQ=\Gal(\ov\QQ/\QQ)$ over a finite field come from modular forms of
prescribed weight, level, and character. This has now been proved by
C.~Khare and J.-P.~Wintenberger; see \cite{KW1,KW2}. Because this result
can be seen as a generalization of Artin Reciprocity to the $\mathrm{GL}_2$
case (over \QQ), we will refer to it as ``Serre Reciprocity.''

Already in \cite{duke}, Serre showed how, given a compatible system of
$\ell$-adic Galois representations and bounds on the weight and level of
the predicted modular forms in characteristic $\ell$, one can use Serre
Reciprocity to obtain results in characteristic zero. We refer to this as
``Serre's method'' and state and prove a generalized form of it in
Section~1 below.

Serre's method allows us to show that certain geometric Galois
representations are modular. Specifically, we show that the representation
obtained from the third \'etale cohomology of a rigid Calabi-Yau threefold
defined over \QQ\ comes from a modular form of weight 4 on $\Gamma_0(N)$.
The proof is an application of Serre's method; it can, in fact, be read off
directly from \cite[Section 4.8]{duke}, which is why one might describe
this short paper as a ``footnote to Serre.'' Recent results allow a
slightly simpler version of the proof.

The observation that the proof of Serre's reciprocity allows us to
establish the modularity of odd irreducible motives of rank two has also
been made independently, in more general terms, by Mark Kisin in
\cite{kisin} (see his Corollary 0.5). 

The question of the modularity of the Galois representations obtained from
Calabi-Yau threefolds over \QQ\ has been much studied, and a large number
of examples are now available; see \cite{Meyer} for a survey. Since current
methods restrict us to low-dimensional Galois representations, many of
the examples involve \text{rigid} Calabi-Yau threefolds, defined below,
simply because in that case the representation is automatically of
dimension two.  

Dieulefait and Manoharmayum have shown in \cite{DiMa} that if $X$ has good
reduction at small primes then it is modular. Richard Taylor showed in
\cite{TaylorCoates} that rigid Calabi-Yau manifolds over \QQ\ are
\emph{potentially modular}, i.e., that there exists a totally real field
$F$ such that the restrictions to $\Gal(\ov\QQ/F)$ of the representations
$\rho_\ell$ are attached to automorphic representations over $F$. These
results were based on the same family of modular lifting theorems that was
used to finally prove Serre reciprocity.

As we will indicate below, the same methods also apply to the non-rigid
case if one can isolate an irreducible two-dimensional ``piece'' of the
cohomology. In this case, we obtain modular forms of weight $4$ and of
weight $2$, in agreement with many examples found by Meyer and others.
In general, however, the middle cohomology groups of non-rigid Calabi-Yau
threefolds do not decompose into products of two-dimensional
pieces.

\vspace{\baselineskip}

\ni\textbf{Acknowledgments}

We would like to thank J.-P. Serre, Richard Taylor, and J.-M. Fontaine for
their willingness to answer questions over email. We thank Matthew Emerton,
Luis Dieulefait, and Matthias Sch\"utt for their very helpful comments on a
preliminary version of this article.  We thank Bert van Geemen
for pointing out the existence of a Galois twist of the Schoen quintic.

The first author thanks Queen's University for its hospitality when this
work was done and Colby College for research funding. The second author
thanks the Natural Sciences and Engineering Research Council of Canada
(NSERC) for their support.

\section{Serre Reciprocity}

Let $G=\Gal(\ov\QQ/\QQ)$ be the absolute Galois group of \QQ, let \FF\ be a
finite field of characteristic $\ell$, and let $\rho: G\to \GL_2(\FF)$ be
an absolutely irreducible representation. We will assume throughout that
$\rho$ is \emph{odd}, that is, we will assume that if $c\in G$ is (any)
complex conjugation, we have $\det\rho(c)=-1$. We will let $S$ be the
finite set of primes such that $\rho$ is unramified at all primes not in
$S$.

In \cite{duke}, Serre associated to any such $\rho$ a triple $(N,k,\eps)$,
where $N$ and $k$ are positive integers, $k\geq 2$, and $\eps$ is a
Dirichlet character modulo $N$. We will briefly recall below how this
triple is obtained, but we refer the reader to \cite{duke} and \cite{Edix}
for details.

In order to avoid technical problems related to finite fields of small
characteristic, we assume $\ell\neq 2,3$. See \cite{Edix} for how to modify
the statement below so that it remains true in those cases.

The result conjectured in \cite{duke} and proved in \cite{KW1,KW2} is:

\begin{theorem}[Serre--Khare--Wintenberger]
  Let $G$, $\ell$, and $\FF$ be as above.  Suppose $\rho:G\to\GL_2(\FF)$ is
  an odd absolutely irreducible representation, and let $(N,k,\eps)$ be the
  Serre parameters attached to $\rho$. Then there exist:
  \begin{itemize}
    \item a cuspidal modular eigenform $f$ on $\Gamma_0(N)$, of weight $k$
      and character $\eps$ and defined over a number field $K$, and 
    \item a prime $\lambda$ of $K$ with residue field \FF
  \end{itemize}
  such that the reduction modulo $\lambda$ of the $\lambda$-adic
  representation attached to $f$ is isomorphic to $\rho$.
\end{theorem}

Part of the power of this result comes from the fact that the triple
$(N,k,\eps)$ is specified in advance in terms of $\rho$, which restricts us
to a finite number of possibilities for the eigenform $f$. It will be
helpful to recall how these parameters are obtained.

The level $N$ is fairly easy to describe: it is the prime-to-$\ell$ part of
the Artin conductor of the representation $\rho$. As such, it is divisible
only by primes $p\in S$, $p\neq\ell$. If we set
\[ N=\prod_{p\in S} p^{e(p)},\] the exponent $e(p)$ is entirely determined
by the image of the inertia group at $p$. In particular, it is useful
to note that if $\rho$ is \emph{tamely} ramified at $p$ then $e(p)=1$.

This choice of the level parameter has as a useful side effect that the
modular form $f$ will necessarily be a \emph{newform}, i.e., it will not
come from a level lower that $N$. (If $f$ did come from a form of lower
level, the $p$-adic representation attached to $f$ would have smaller
conductor, and therefore so would the mod $p$ representation.)

The weight $k$ is the most delicate of the three parameters. It depends
only on the image of the inertia group at $\ell$, but the recipe for
computing $k$ is complicated; see \cite{duke} and \cite{Edix} for details.
We will use Serre's normalization of the weight, so that $k\geq 2$.

Finally, the character $\eps$ is determined by the formula
\[ \det\rho = \eps\chi_{\ell}^{k-1},\]
where $\chi_\ell$ is the (reduction mod $\ell$ of the) $\ell$-adic
cyclotomic character. Notice that this formula determines $k$ modulo
$\ell-1$. 

While Serre Reciprocity is a statement about Galois representations over
finite fields, it can often be used to show the modularity of
representations in characteristic zero as well. The idea is due to Serre
himself; it will perhaps be useful to have a formalized version of it.

\begin{theorem}[Serre's Method]\label{method}
Fix a number field $K$, and let $\lambda$ run over primes of $K$. For each
$\lambda$, let $K_\lambda$ denote the completion of $K$ at $\lambda$ and
let $\kappa(\lambda)$ be the residue field. Let $\ell$ be the
characteristic of $\kappa(\lambda)$. 

Fix a finite set $S$ of primes in \QQ. For each $p\notin S$, let
$\mathrm{Frob}_p$ be a choice of arithmetic Frobenius element in
$\Gal(\ov\QQ/\QQ)$.

Suppose we have, for each $\lambda$, a two-dimensional $K_\lambda$ vector
space $V_\lambda$ with a continuous action of $G=\Gal(\ov\QQ/\QQ)$. This
gives a continuous representation $\rho_\lambda$. Assume $\rho_\lambda$ is
odd and unramified outside $S\cup \{\ell\}$.

For each $\lambda$ we can find a $G$-stable lattice, reduce modulo
$\lambda$, and semi-simplify if necessary to obtain a semisimple odd
two-dimensional Galois representation
\[ \ov\rho_\lambda: G \longrightarrow \GL_2(\kappa(\lambda))\]
unramified outside $S\cup \{\ell\}$.

Fix an infinite set $I$ of primes in $K$. Suppose we can show that:
\begin{enumerate}
\item For all $\lambda\in I$, the representation
  $\ov\rho_\lambda$ is absolutely irreducible.
\item There exists a family $Q_p(X)= X^2-A_pX+D_p \in K[X]$ of polynomials
  of degree $2$, indexed by primes $p\in\QQ$, $p\notin S$, such that for all
  $\lambda\in I$ and (given $\lambda$) all $p\notin S\cup \{\ell\}$, the
  characteristic polynomial of $\mathrm{Frob}_p$ acting on $V_\lambda$ is
  equal to $Q_p(X)$.
\item There exists an integer $k_0$ such that for all $\lambda \in I$ the
  Serre weight $k_\lambda$ attached to $\ov\rho_\lambda$ satisfies
  $1<k_\lambda\leq k_0$.
\item There exists an integer $N_0$ such that for all $\lambda \in I$ the
  Serre level $N_\lambda$ attached to $\ov\rho_\lambda$ is a divisor of
  $N_0$. We choose $N_0$ to be minimal with this property.
\end{enumerate}

Then there exists a cuspidal Hecke eigenform form $f$ (new of level $N$
dividing $N_0$, weight $k$ less than or equal to $k_0$, defined over $K$)
such that for all $\lambda$ the $\lambda$-adic representation
$\rho_{f,\lambda}$ attached to $f$ is isomorphic to $\rho_\lambda$.
\end{theorem}

\begin{proof}
  
  We may, and will, assume that the set $S$ has been chosen to be as small
  as possible, so that for every $p\in S$ there is at least one $\lambda$
  that does not divide $p$ and such that $\rho_\lambda$ is ramified at $p$.
  Note also that since $\ov\rho_\lambda$ is absolutely irreducible for some
  $\lambda$, so is $\rho_\lambda$. (Because characters can be inserted into
  compatible systems, it follows from the compatibility condition (2) that
  $\rho_\lambda$ is absolutely irreducible for \emph{all} $\lambda$.  We
  will not actually use this, and once we have shown that the
  representation comes from a cuspform, it will follow that all the
  $\rho_\lambda$ are absolutely irreducible.)

  Choose $\lambda\in I$ and apply Serre Reciprocity to $\ov\rho_\lambda$.
  We get an eigenform $f_\lambda$ of weight $k_\lambda$ less than or equal
  to $k_0$, character $\eps$, and level dividing $N_0$. A priori, $f$ may
  be defined over an extension of $K$ whose residue field at a prime
  $\lambda'$ over $\lambda$ is still $\kappa(\lambda)$. The fact that
  $f_\lambda$ corresponds to $\rho_\lambda$ tells us that
  \[ A_p \equiv a_p(f_\lambda) (mod \,\,\lambda') \]
  and
  \[ D_p \equiv \eps(p)p^{k_\lambda-1}  (mod \,\,\lambda')\]
  for all $p\notin S\cup\{\ell\}$.
  
  Since the set of all eigenforms of weight bounded by $k_0$ and level
  dividing $N_0$ is \emph{finite} and there are infinitely many $\lambda\in
  I$, there must exist a modular form $f$ such that $f_\lambda=f$ for
  infinitely many $\lambda$. Let $k$ be the weight of $f$. But then, for
  each $p\notin S$ we will have
  \[ A_p \equiv a_p(f) (mod \,\,\lambda) \]
  and
  \[ D_p \equiv \eps(p)p^{k-1}  (mod \,\,\lambda)\]
  for infinitely many $\lambda$.
  
  This implies that in fact $A_p=a_p$ and $D_p=\eps(p)p^k$ for all $p\notin
  S$. Since we know $f$ is a newform, this is enough to show that $f$ is
  the eigenform we wanted to find and (together with the minimality of $S$)
  implies, in particular, that it has coefficients in $K$.
\end{proof}

In the case of representations coming from geometry, the representations
will typically be obtained from the (dual of the) \'etale cohomology of an
algebraic variety $X$ defined over $\QQ$. The field $K$ is then just $\QQ$.
The set $S$ is then contained in the set of primes of bad reduction for $X$
and the existence of the $Q_p(X)$ follows from the Weil Conjectures as
proved by Deligne.

For an example with $K$ a totally real field, consider the case of abelian
varieties with real multiplication; see \cite[Section 4.7]{duke}

\section{Modularity of rigid Calabi-Yau threefolds over $\QQ$}

We want to apply Serre's method to the representation obtained from
the middle \'etale cohomology of a rigid Calabi-Yau threefold defined
over \QQ. We recall the definitions.

\begin{definition} \label{CYmfd}
Let $X$ be a smooth projective threefold defined over $\CC$.
We call $X$ a \emph{Calabi-Yau threefold} if
\begin{enumerate}
\item $H^1(X,\mathcal{O}_X)=H^2(X,\mathcal{O}_X)=0$, and
\item $K_X:=\wedge^3 \Omega_X^1\simeq \mathcal {O}_X$, that is, the
  canonical bundle is trivial.  
\end{enumerate}
\end{definition}

As usual, we define the Hodge numbers
\[ h^{i,j}(X):=\mathrm{dim}_{\CC} H^j(X,\Omega^i_X).\]
By complex conjugation, $h^{i,j}(X)=h^{j,i}(X)$, and by Serre duality,
$h^{i,j}(X)=h^{3-j,3-i}(X)$ for $0\leq i,\, j\leq 3$. 
The Hodge decomposition gives 
\[ h^k(X)=\mathrm{dim}_{\CC}H^k(X,\CC)=\sum_{i+j=k} h^{i,j}(X).\] 
The number $h^k(X)$ is called the $k$-th Betti number of $X$ and often
denoted $B_k(X)$. 

If $X$ is Calabi-Yau, then the first condition implies that
\[ h^{1,0}(X)=h^{2,0}(X)=0,\]
and the second condition, together with Serre duality, yields
\[ h^{3,0}=h^{0,3}=1.\] We can summarize all this by drawing the ``Hodge
diamond'' of $X$:

\[
\begin{array}{c@{}c@{}c@{}c@{}c@{}c@{}c@{}cl}
  &   &            &          1 &            &   &   & \qquad\qquad & h^0(X) = 1\\
  &   & 0          &            &          0 &   &   & & h^1(X) = 0\\
  & 0 &            & h^{1,1}(X) &            & 0 &   & & h^2(X) = h^{1,1}(X)\\
1\phantom{123} &   & h^{2,1}(X) &            & h^{1,2}(X) &   & \phantom{123}1 & & h^3(X) = 2(1+h^{2,1}(X))\\
  & 0 &            & h^{2,2}(X) &            & 0 &   & & h^4(X) = h^{2,2}(X) = h^{1,1}(X)\\
  &   & 0          &            &          0 &   &   & & h^5(X) = 0\\
  &   &            &          1 &            &   &   & & h^6(X) = 1
\end{array}
\]

Calabi-Yau threefolds are K\"ahler manifolds, so $h^{1,1}(X)>0$.  All
$2$-cycles on Calabi-Yau threefolds are algebraic, as follows from the
Lefschetz $(1,1)$ theorem that $H^2(X,\ZZ)\cong \mathrm{Pic}(\ov X).$ In
particular, $h^{1,1}(X)=h^2(X)=\mathrm{rk}\,\mathrm{Pic}(\ov X)$.

\begin{definition} \label{rigid}
Let $X$ be a Calabi-Yau threefold defined over $\CC$. We say
that $X$ is \emph{rigid} if $h^{2,1}(X)=h^{1,2}(X)=0$, so that
$h^3(X)=2$.
\end{definition}

The name ``rigid'' comes from the fact that the space of deformations of a
Calabi-Yau manifold has dimension $h^{1,2}(X)$.  There are more than $50$
known examples of rigid Calabi-Yau threefolds (up to birational equivalence
over \CC). It is still an open problem to decide whether the number of such
examples is finite up to birational transformation over \CC.

Since we are interested in Galois representations, we focus on Calabi-Yau
threefolds defined over $\QQ$. Of course, a given Calabi-Yau threefold over
\CC\ may well have many different realizations over \QQ. Notice that the
fact that $X$ is a rigid Calabi-Yau manifold is independent of the choice
of model over \QQ, but the Galois representation (and therefore the modular
forms we will find) depend strongly on that choice. We will comment further
on this below.

Let $X$ be a rigid Calabi-Yau threefold defined over \QQ. Then $X$ always
has a model defined over $\ZZ$; we assume one has been chosen and fixed. We
will apply Theorem~\ref{method} with $K=\QQ$ and $S$ the (finite) set of
primes at which $X$ has bad reduction.  Let $\ov X=X\tensor\ov{\QQ}$.  For
each prime $\ell$ in \QQ, let
\[ V_\ell = H^3(\ov X, \QQ_\ell)\spcheck. \]
(We need to dualize because we want to work with the arithmetic Frobenius.)
We know that if $p\notin S\cup \ell$, this representation will be
unramified at $p$.

The assumption that $X$ is rigid means that that $V_\ell$ is two
dimensional and that its Hodge decomposition is of the form $(3,0)+(0,3)$.
By Pontryagin duality, we have
\[ \det\rho_\ell = \chi_\ell^3\]
so that $\rho_\ell$ is odd. Let $\ov\rho_\ell$ be the representation
obtained by reducing modulo $\ell$.

In \cite[Section 4.8]{duke}, Serre checked that conditions (1) and (2)
above hold for sufficiently large $\ell$. A theorem of Fontaine (see also
\cite{Edix}) shows that for all large enough $\ell$ the Serre weight
parameter will be $k=4$.  

In order to verify the condition on the level, Serre used a bound for the
Artin conductor proved in \cite[Section 4.9]{duke}: under certain
congruence conditions on $\ell$, the conductor $N$ is a divisor of
\[ N_0 = \prod_{p\in S} p^{e(p)},\]
where $e(2)=8$, $e(3)=5$, and $e(p)=2$ for all other primes $p\in S$. We
can therefore let $I$ be the (infinite) set of primes $\ell$ that satisfy
Serre's congruence conditions.

This can now be simplified by using the results in \cite{TaylorCoates}.
Since the Hodge numbers are $0$ and $3$ and we know that the representation
is crystalline at all $p\notin S$ (because $X$ has good reduction at all
such primes), the $\rho_\ell$ form what Taylor calls a \emph{weakly
  compatible} system of representations; by Theorem~A in
\cite{TaylorCoates}, the system must in fact be \emph{strongly} compatible,
which implies that the conductor $N$ is \emph{independent} of $\ell\notin
S$. (We thank Luis Dieulefait for pointing this out to us.) Hence we can
take our infinite set to be all primes not in $S$.

Theorem~\ref{method} then gives our result: 

\begin{theorem}
  Let $X$ be a rigid Calabi-Yau threefold defined over \QQ, and use the
  notations above. Then there exists a Hecke eigenform $f$ of weight $4$,
  level dividing $N$, and trivial character such that $\rho_\ell$ is
  equivalent to $\rho_{f,\ell}$ for all $\ell$.
\end{theorem}

In other words, all rigid Calabi-Yau threefolds defined over \QQ\ are
modular. In particular, this implies that the $L$-function corresponding to
the third \'etale cohomology of such a threefold is the same as that of a
modular form of weight $4$, and hence is holomorphic and satisfies a
functional equation relating values at $s$ to values at $4-s$.

Notice that while we do not need to use Serre's bound on the level for the
argument, a posteriori the bound will apply to the level of the form $f$.
This is in fact the bound obtained by Dieulefait in \cite{Di}.

Serre's method is applicable, as he shows in \cite{duke}, to all
odd-dimensional smooth algebraic varieties whose middle-dimensional
cohomology is of dimension two and of Hodge type $(*,0)+(0,*)$.

\section{The Non-Rigid Case}

The reason to focus on the rigid case is, of course, that we get a Galois
representation of dimension two, which should then come from a modular
form. Higher dimensional representations should be automorphic, but the
type of corresponding automorphic representation we expect to find will
depend on that dimension.

If we drop the assumption that the Calabi-Yau manifold $X$ is rigid, then
$h^3$ will not be equal to two. It is still possible, nevertheless, that
the Galois representation attached to the third cohomology contains an
irreducible subrepresentation of dimension two. If such a subrepresentation
occurs in $H^3(\ov X,\QQ_\ell)$ for every $\ell$ and the resulting Galois
representations are (weakly) compatible, the same argument will apply. Such
a system of compatible representations is usually described as a submotive
of rank two.

If the submotive happens to be the $(3,0)+(0,3)$ part, exactly the same
argument will show that it is modular, i.e., the subrepresentation of
dimension two over $\QQ_\ell$ will be isomorphic to the $\ell$-adic
representation attached to a modular form $f$ of weight $4$.

If the submotive $M_\ell$ occurs instead in the $(2,1)+(1,2)$ part, the
method described above is not directly applicable, because the Serre weight
$k_\ell$ will in general be $\ell+3$, and hence not bounded. This can be
easily fixed, however, by twisting: $M_\ell\tensor\chi_\ell^{-1}$ has Hodge
numbers $(1,0)+(0,1)$, and the argument above will show that it corresponds
to a modular form of weight $2$. Hence $M_\ell=V_\ell(f)\tensor\chi_\ell$
is a Tate twist of the representation coming from such a form of weight
two.

There are several (proved and conjectural) examples of this in
\cite{Meyer}. Many of them are of the type studied in \cite{HV}, namely,
Calabi-Yau threefolds containing a large number of elliptic ruled surfaces.
To be specific, let
\[ V_\ell = H^3(\ov X,\QQ_\ell)\spcheck.\]
The examples in \cite{HV} and \cite{Meyer} look like
\[ V_\ell \cong V_\ell(f) \oplus \left[V_\ell(g_1)\otimes\chi_\ell\right]
\oplus \left[V_\ell(g_2)\otimes\chi_\ell\right] \oplus \dots \oplus
\left[V_\ell(g_k)\otimes\chi_\ell\right],\] where $h^3(X)=2+2k$, $f$ is a
modular form of weight $4$, the $g_i$ are all modular forms of weight $2$,
and $V_\ell(h)$ is the $\ell$-adic representation attached to a modular
eigenform $h$. In many cases, Meyer finds examples where the $g_i$ are in
fact all the same; see, for example, pages 23--24 of \cite{Meyer}.  Another
example can be found in \cite{Lee}, where we have $h^3=4$ and the
representation splits into two two-dimensional components.

Finally, if the decomposition of the cohomology representation takes place
only after extending scalars, but we still obtain weakly compatible
systems, Serre's method applies to deduce modularity (after twisting) of
all of them, the only difference being that the relevant modular forms are
no longer defined over \QQ. We thank the referee for pointing that out to
us. 

\section{Some speculations}

Let $X$ be a rigid Calabi-Yau threefold defined over \QQ, and let $f$ be
the associated modular form of weight $4$.

1) The level $N$ of the form $f$ is going to be a delicate arithmetic
invariant of $X$ (over \QQ, rather than over the algebraic closure). The
primes dividing $N$ should be primes at which $X$ has bad reduction in
every model of $X$ over \ZZ, but it is unclear whether $N$ \emph{must} be
divisible by all such primes. In addition, the precise power of such primes
that occurs in $N$ presumably depends on the type of singularities, but we
do not know how that should work.

2) Suppose we have $X$ and its modular form $f$ of level $N$. Then if we
twist $f$ by a quadratic character of level $d$, we get another eigenform
of weight $4$ and level dividing $Nd^2$. Will this form be attached to
another rigid Calabi-Yau threefold over \QQ? If so, then there must exist
an algebraic correspondence between $X$ and $X_d$ that is defined over
$\QQ(\sqrt{d})$. We might even hope that $X_d$ is a Galois twist of $X$, so
that $X\tensor\QQ(\sqrt{d}) \cong X_d\tensor\QQ(\sqrt{d})$. (Note that the
existence of $X_d$ which are Galois twists is compatible with the
conjecture that there are only finitely many rigid Calabi-Yau threefolds
over \CC{} up to birational equivalence.)

Helena Verrill has found an example where such a Galois twist $X_d$ can be
constructed (see \cite{Y03}); does such a ``twist'' always exist? In
\cite{Meyer}, Meyer conjectured that the answer is yes. An interesting test
case is the Schoen quintic
\[ X_0^5 +X_1^5 +X_2^5 +X_3^5 +X_4^5 = 5X_0X_1X_2X_3X_4.\]
This is a singular threefold, and resolving those singularities produces a
rigid Calabi-Yau threefold $X$ that is known (see \cite{Meyer} and the
references therein) to be associated to a modular form of weight $4$ and
level $25$. Can one contruct the requisite Galois twists?

Bert van Geemen has informed us that the answer is
``yes''.  Since the Schoen quintic $X$ has an
automorphism $\phi$ of order $2$ defined over $\QQ$, 
which acts on $H^3(X)$ by $-1$, and we can use
this to twist the quintic $X$.
Let $K=\QQ(\sqrt{d})$ and let $X_K$ be the quintic defined over $K$.
Descend $X_K$ back to $\QQ$ by taking the quotient by the
automorphism which is $\phi$ on the coordinates and which is
the non-trivial automorphism of $K$ on the scalars. Then
we get a quintic $X_d$ on whose $H^3(\overline{X_d},\QQ_{\ell})$ 
the Galois representation is the twist of the one on 
$H^3(\overline{X},\QQ_{\ell})$.

In fact, the automorphism $\phi$ of order $2$ is given, for instance, 
explicitly by 
$$\phi(X_0)=X_1,\, \phi(X_1)=X_0,\, \phi(X_i)=X_i\quad\mbox
{for $i=2,3,4$.}$$
Put $U=X_0+X_1$ and $V=X_0=X_1$. Then the equation for
the quintic equation can be written as a polynomial in $U$ and $V^2$
as follows:
$$U^5+10U^3V^2+5UV^4+16(X_2^5+X_3^5+X_4^5)-20(U^2-V^2)X_2X_3X_4=0.$$
Now replace $V$ by $\sqrt{d}V$, then we obtain the quintic
equation for $X_d$:   
$$U^5+10dU^3V^4+16(X_2^5+X_3^5+X_4^5)-20(U^2-dV^2)X_2X_3X_4=0.$$
By counting points on $X_d$, we can see explicitly that the
Galois representation has been twisted.

3) Can we reverse this process? In other words, given an eigenform $f$
of weight $4$ on $\Gamma_0(N)$ and defined over \QQ, does there exist
a rigid Calabi-Yau threefold $X$ corresponding to $f$? Since Barry
Mazur first called attention to this question, it is known as
\textit{Mazur's problem}. 

Of course, Mazur's problem is also connected to the issue of whether there
are infinitely many different birational equivalence classes of rigid
Calabi-Yau threefolds. If the answer to Mazur's question is ``yes,'' then
we can translate the question to the setting of modular forms, where it
becomes the question of understanding whether there are, up to twists,
infinitely many modular forms of weight $4$ (of any level) that are defined
over \QQ.

\end{document}